\documentclass[10pt]{elsarticle}
\usepackage[cp1251]{inputenc}
\usepackage{amsmath}
\usepackage{amssymb}
\usepackage{amsfonts}
\usepackage{mathtools}
\usepackage{dsfont}
\usepackage{graphicx}
\usepackage{floatflt,epsfig}
\usepackage{lineno,hyperref}
\usepackage{enumerate}
\usepackage{colortbl}
\usepackage{array,tabularx,tabulary,booktabs}
\usepackage{longtable}
\usepackage{multirow}
\usepackage{wrapfig}
\usepackage{subcaption}
\usepackage{caption}
\usepackage[table]{xcolor}
\usepackage{adjustbox}
\usepackage{rotating}
\usepackage{arydshln}

\usepackage{tabularray}
\newcolumntype{^}{>{\currentrowstyle}}

\journal{arXiv}
\newtheorem{lemma}{Lemma}

\newtheorem{theorem}{Theorem}

\newcommand{\proof}{\medskip\noindent{\bf Proof.~}}
\begin{document}
\renewcommand{\abstractname}{Abstract}
\renewcommand{\refname}{References}
\renewcommand{\arraystretch}{0.9}
\thispagestyle{empty}
\sloppy

\begin{frontmatter}
\title{Distinct eigenvalues of the Transposition graph}

\author[01,02,03]{Elena~V.~Konstantinova}
\ead{e\_konsta@math.nsc.ru}

\author[01,02]{Artem Kravchuk}
\ead{artemkravchuk13@gmail.com}

\address[01]{Sobolev Institute of Mathematics, Ak. Koptyug av. 4, Novosibirsk 630090, Russia}
\address[02]{Novosibirsk State University, Pirogova str. 2, Novosibirsk, 630090, Russia}
\address[03] {Three Gorges Mathematical Research Center, China Three Gorges University, 8 University Avenue, Yichang 443002, Hubei Province, China}

\begin{abstract}
Transposition graph $T_n$ is defined as a Cayley graph over the symmetric group generated by all transpositions. It is known that all eigenvalues of $T_n$ are integers. Moreover, zero is its eigenvalue for any $n\geqslant 4$. But the exact distribution of the spectrum of the graph $T_n$ is unknown. In this paper we prove that integers from the interval $[-\frac{n-4}{2}, \frac{n-4}{2}]$ lie in the spectrum of $T_n$ if $n \geqslant 19$.
\end{abstract}

\begin{keyword}
Transposition graph; integral graph; spectrum; 
\vspace{\baselineskip}
\MSC[2010] 05C25\sep 05E10\sep 05E15
\end{keyword}
\end{frontmatter}

\section{Introduction}\label{sec1}

The {\it Transposition graph} $T_n$ is defined as a Cayley graph over the symmetric group $\mathrm{Sym}_n$ generated by all transpositions. The graph $T_n, n\geqslant 2$, is a connected bipartite $\binom{n}{2}$-regular graph of order $n!$ and diameter $(n-1)$~\cite{K08}. It is an edge--transitive graph but not distance--regular, and hence not distance--transitive graph. It was shown in~\cite{KL20} that the Transposition graph is integral which means that all eigenvalues of its adjacency matrix are integers~\cite{HS74}. Since $T_n$ is bipartite then its spectrum $Spec(T_n)$ is symmetric with respect to zero, where the spectrum of a graph is defined as a multiset of distinct eigenvalues together with their multiplicities~\cite{BH12}. 

Let $G$ be a finite group with an identity element $1_G$, and $S$ be its generating set. Then a Cayley graph $\Gamma=Cay(G,S)$ is called {\it normal} if its generating set $S$ is closed under conjugation, i.~e. $S$ is the union of conjugacy classes of $G$.  The following theorem enables to compute eigenvalues for any normal Cayley graph $\Gamma$ in terms of complex character values of $G$.  

\begin{theorem}\label{Z-88} {\rm \cite[Theorem~1]{Z88}} Let $G$ be a finite group with $s$ conjugacy classes and let $\{\chi_1,\chi_2,\ldots \chi_s \}$ be the set of all irreducible complex characters of $G$. Then the eigenvalues $\lambda_i, \ i=1,2,\ldots,s$, of any normal Cayley graph $\Gamma=Cay(G,S)$ are given by the following expression:

\begin{equation}\label{eigen_expr}
\lambda_i=\sum_{g\in S} \frac{\chi_i(g)}{\chi_i(1_G)},  
\end{equation}
\end{theorem}

It was shown in~\cite{KY97} that the conditions of Theorem~\ref{Z-88} hold for the Transposition graph $T_n=Cay(\mathrm{Sym}_n,T)$, where $T$ is the set of all transpositions. Moreover, the following useful expressions were obtained.

It is well-known fact (see~\cite{Sa01}) that there is one-to-one correspondence between the irreducible complex characters $\chi_i(g)$ and $\chi_i(1_G)$, $i=1,2,\ldots,p(n)$, of the symmetric group $\mathrm{Sym}_n$ and the partitions of $n$, where $p(n)$ is the number of partitions of $n$. Let a nonincreasing sequence $(n_1,n_2,\ldots,n_k), \ k\geqslant 1$, where $\sum_{j=1}^k n_j=n$, be the integer partition ${\bf i}=(n_1,\dots,n_k)\vdash n$ of $n$ corresponding to an irreducible complex character $\chi_i$. Then the following expression holds:
\begin{equation}\label{Formula3-KY}
\frac{\chi_i(\tau)}{\chi_i(I_n)}=\sum_{j=1}^k \frac{n_j(n_j-2j+1)}{n(n-1)}, 
\end{equation}
where $\tau$ is a transposition and $I_n$ is the identity permutation. Since the generating set $T$ of $T_n$ consists of $\frac{n(n-1)}{2}$ transpositions then equations~(\ref{eigen_expr}) and~(\ref{Formula3-KY}) give an expression for an eigenvalue $\lambda_{\bf i}$ corresponding to the partition ${\bf i}$:

\begin{equation}\label{transp_eigen}
    \lambda_{\bf i} = \sum_{j=1}^k \frac{n_j(n_j-2j+1)}{2}. 
\end{equation}

This formula for the eigenvalue corresponding to ${\bf i}$ is precisely the same as the normalized character on the sum of all transpositions found in \cite[Equation~(4.1)]{Cesi2010} and is a simplified version of the character of a transposition found in \cite[Lemma~7]{DS81}.

The following result about the spectrum of the graph $T_n$ is known.

\begin{theorem} \label{KK22} {\rm ~\cite[Theorem~2]{KK22}}
 For any integer $k\geqslant 0$, there exists $n(k)$ such that for any $n\geqslant n(k)$ and any $m \in \{0, \dots, k\}$, $m \in Spec(T_n)$. 
 \end{theorem}

In the proof of this theorem $n(k)$ is chosen as $n(k)=10k+4$. It follows that all integers from the interval $[-\frac{n-4}{10}, \frac{n-4}{10}]$ lie in the spectrum of $T_n$.

In this paper we improve the above result by the following theorem. 

\begin{theorem}\label{segment_0_n2_theorem}
For  $n \geqslant 19$, all integers from the segment $[-\frac{n-4}{2}, \frac{n-4}{2}]$ lie in the spectrum of $T_n$.
\end{theorem}

The paper is organized as follows. In Section~\ref{sec2} we give the proof of Theorem~\ref{segment_0_n2_theorem}. This proof relies on a few technical lemmas proved in Section~\ref{tech-lemmas}. The summaries of the lemmas are collected in Table~\ref{main-table}. These lemmas mainly use~(\ref{transp_eigen}) to show that the given partitions correspond to the desired eigenvalues.

\section{Proof of Theorem~\ref{segment_0_n2_theorem}}\label{sec2}

Since $T_n$ is bipartite it is sufficient to prove that all integers from the segment $[0, \frac{n-4}{2}]$ lie in the spectrum of $T_n$. To prove this statement, we divide the segment $[0, \frac{n-4}{2}]$ into subsegments $S = \bigcup_{i}s_i$, and specific integer numbers $A=\{a_0, \dots, a_k \}$ depending on the parity of $n$, where $S \cup A$ contains all integers of the interval $[0, \frac{n-4}{2}]$, as follows:

\setlength{\arrayrulewidth}{0.1mm}
\setlength{\tabcolsep}{8pt}
\renewcommand{\arraystretch}{1.5}
\begin{center}
\begin{tabular}{ |c|c|c| } 
 \hline
 $n$ & $S$ & $A$ \\ 
 \hline
 $odd$ & $\{[1, \frac{n-3}{4}], [\frac{n+3}{4}, \frac{n-1}{2}]\}$ & $\{0, [\frac{n-3}{4}] + 1\}$ \\
 \hline
 $even$ & $\{[2, \frac{n-4}{4}], [\frac{n+2}{2}, \frac{n-4}{2}]\}$ & $\{0, 1, [\frac{n-3}{4}] + 1\}$ \\ 
 \hline

\end{tabular}
\end{center}

Next, we show that all integers from each $s\in S$ lie in the spectrum of $T_n$. For each $s$, we find a family of partitions for which by~(\ref{transp_eigen}) it is shown that all the integers from $s$ are elements of the spectrum of $T_n$. For each integer number $a \in A$, we find a particular partition for which by~(\ref{transp_eigen}) the eigenvalue $a \in Spec(T_n)$ is obtained. Thus, we show that all integers from the segment $[0,\frac{n-4}{2}]$ lie in the spectrum of $T_n$, and since $T_n$ is an integral graph, there are no other eigenvalues in the interval $[0, \frac{n-4}{2}]$.

To realise the strategy above, we prove a few technical lemmas. Summary of these lemmas are presented in Table~\ref{main-table} which lists $s \in S$ and $a \in A$ in the first column. Whenever a segment is listed this means that all integers from the segment are in the spectrum of $T_n$. The second column lists the specific partitions, or families of partitions, which correspond to a given number or a given segment. The notation $(n_1, \dots, k\times t, \dots, n_k)$ means that number $k$ is repeated $t$ times. The specific eigenvalues (e.g., $0$) correspond to the specific partitions, and the segments (e.g., $[1, \frac{n-3}{4}]$) correspond to the families of partitions. If we choose an arbitrary integer value $\lambda$ from the segment from the first column and substitute it in the expression for the family of partitions from the second column we obtain a partition corresponding to the eigenvalue $\lambda$ of $T_n$. It is worth noting that some partitions may depend not only on the parity of $n$ but also on the remainder of $n$ divided by 4.

For each family of partitions, we prove a specific lemma to show that this family gives the corresponding eigenvalues. A reference to this proof is given in the corresponding column of Table~\ref{main-table}. Each lemma has its own restrictions on the number $n$. These constraints are given in the last column of Table~\ref{main-table}. 

Each line of Table~\ref{main-table} can be taken as the result of a single technical lemma. Some of the lemmas were proved in~\cite{KK22}. The remaining lemmas are proved in Section~\ref{tech-lemmas}. Since these proofs are almost identical and a little bit cumbersome they are moved to a separate section so that the reader can refer to them if necessary but this would not lead away from the main line of the proof. 

Thus,  all that remains is to choose the minimum value $n$ for which all the lemmas hold. It is easy to see from Table~\ref{main-table} that all lemmas hold for $n \geqslant 19$ since the maximum number for which one of the lemmas does not hold is $n=18$, namely, Lemma~\ref{even_mod2}. 

\hfill $\square$

\setlength{\arrayrulewidth}{0.1mm}
\setlength{\tabcolsep}{8pt}
\renewcommand{\arraystretch}{2}

\begin{sidewaystable}[!h]
    \centering
   \begin{tabular}{ |c|c|c|c| } 
\hline
Eigenvalues, $\lambda$ & Partitions & Limitations & Proof  \\
\hline
\multirow{2}{4em}{$\ \ \ \ \ 0$} & $\left(\frac{n+1}{2}, 1\times\frac{n-1}{2}\right)$ &$n$ is odd, $n\geqslant 1$ & \multirow{2}{4em}{\rm ~\cite{KK22}, Lemma~3}\\ 
& $\left(\frac{n}{2}, 2, 1\times \frac{n-4}{2}\right)$ & $n$ is even, $n\geqslant 4$ & \\
\hline
$[1, \frac{n-3}{4}]$ & $(\frac{n-2\lambda+1}{2}, \lambda+2, 2\times (\lambda-1), 1\times \frac{n-4\lambda-1}{2})$ & $n$ is odd, $n\geqslant 7$ & \rm ~\cite{KK22}, Lemma~5\\
\hdashline
$1$ & $\left(\frac{n-6}{2},4,4,2,1\times \frac{n-14}{2}\right)$ & $n$ is even, $n \geqslant 14$ & \rm ~\cite{KK22}, Lemma~4 \\

$[2, \frac{n-4}{4}]$ & $(\frac{n-2\lambda}{2}, \lambda+2, 3, 2\times (\lambda-2), 1\times \frac{n-4\lambda-2}{2})$ & $n$ is even, $n\geqslant 12$ & Lemma~\ref{even_0_n-44} \\

\hline
$[\frac{n+3}{4}, \frac{n-1}{2}]$ &  $(\lambda + 1, \frac{n + 3 - 2\lambda}{2}, 2\times (\frac{n-1}{2} - \lambda), 1\times (\frac{4\lambda - n - 3}{2})$ & $n$ is odd, $n \geqslant 5$ & Lemma~\ref{odd_n2}\\ 
$[\frac{n+2}{2}, \frac{n-4}{2}]$& $(\lambda + 1, \frac{n + 2 - 2\lambda}{2}, 3,  2\times (\frac{n-4}{2} - \lambda), 1\times (\frac{4\lambda - n - 2}{2})$ & $n$ is even, $n \geqslant 10$ & Lemma~\ref{even_n2}\\ 
\hline

\multirow{4}{4em}{$ [\frac{n-3}{4}] + 1$} & 
     $(\frac{n}{4}, \frac{n}{4}, 5, 4, 2 \times (\frac{n-20}{4}), 1)$ & $n\equiv 0 \ (mod \ 4), n\geqslant 20$ & Lemma~\ref{even_mod0} \\
    &  $(\frac{n + 3}{4}, \frac{n + 3}{4}, 4,  2\times (\frac{n - 13}{4}), 1)$ & $n\equiv 1 \ (mod \ 4), n \geqslant 13$ & Lemma~\ref{odd_mod1} \\
    & $(\frac{n + 2}{4}, \frac{n - 2}{4}, 5, 4, 2 \times (\frac{n-22}{4}), 1, 1)$ & $n\equiv 2 \ (mod \ 4), n\geqslant 22$ & Lemma~\ref{even_mod2} \\
    & $(\frac{n + 1}{4}, \frac{n + 1}{4}, 5, 3, 2\times (\frac{n - 19}{4}), 1)$ & $n\equiv 3 \ (mod \ 4), n\geqslant 19$ & Lemma~\ref{odd_mod3} \\
    \hline
\end{tabular}  
    \caption{ \label{main-table}  Summary of the technical lemmas.}
\end{sidewaystable}

\section{Technical lemmas}\label{tech-lemmas}

\begin{lemma}\label{kk_22_l5} {\rm ~\cite[Lemma~5]{KK22}}
If $n \geqslant 7$ is odd then the partition $(\frac{n-2\lambda+1}{2}, \lambda+2, 2\times (\lambda-1), 1\times \frac{n-4\lambda-1}{2})$ corresponds to the eigenvalue $\lambda \in \mathbb{N}$, where $1\leqslant \lambda \leqslant \frac{n-3}{4}$.
\end{lemma}

\begin{lemma}\label{eigens_01}{\rm ~\cite[Corollary of Lemma~3 and Lemma~4]{KK22}}
If $n \geqslant 14$ then the eigenvalues zero and one lies in the spectrum of $T_n$.
\end{lemma}

\begin{lemma}\label{replace2_3}
If a partition $(n_1, n_2, 2, n_4, \dots, n_k), n_2 \geqslant 3$, corresponds to an eigenvalue $\lambda$ then the partition $(n_1, n_2, 3, n_4,  \dots, n_k)$ corresponds to the eigenvalue $\lambda$ too.
\end{lemma}

\proof The proof follows directly from~(\ref{transp_eigen}).

\hfill $\square$

\begin{lemma}\label{even_0_n-44}
If $n\geqslant 12$ is even then the partition $(\frac{n-2\lambda}{2}, \lambda+2, 3, 2\times (\lambda-2), 1\times \frac{n-4\lambda-2}{2})$ corresponds to the eigenvalue $\lambda \in \mathbb{N}$, where $2 \leqslant \lambda \leqslant \frac{n - 4}{4}$.
\end{lemma}

\proof If $\lambda \geqslant 2$ then for each partition of odd $n$ from Lemma~\ref{kk_22_l5} let us replace two by three at the third position. By Lemma [\ref{replace2_3}, this gives the partition of $n+1$ which corresponds to the eigenvalue $\lambda$. It is easy to see that this results in a partition $(\frac{n-2\lambda}{2}, \lambda+2, 3, 2\times (\lambda-2), 1\times \frac{n-4\lambda-2}{2})$ for even $n$. This partition holds when  $2 \leqslant \lambda \leqslant \frac{n - 4}{4}$, and the last inequality is valid for $n \geqslant 12$.

\hfill $\square$

\begin{lemma}\label{odd_n2}
    If $n\geqslant 5$ and $n$ is odd then the partition $(\lambda + 1, \frac{n + 3 - 2\lambda}{2}, 2\times (\frac{n-1}{2} - \lambda), 1\times (\frac{4\lambda - n - 3}{2}))$ corresponds to the eigenvalue $\lambda \in \mathbb{N}$, where $\frac{n + 3}{4} \leqslant \lambda \leqslant \frac{n - 1}{2}$.
\end{lemma}

\proof By~(\ref{transp_eigen}) we have:
$$\lambda_{\left(\lambda + 1, \frac{n + 3 - 2 \lambda}{2}, 2\times(\frac{n-1}{2} - \lambda), 1\times (\frac{4\lambda - n - 3}{2}) \right)}=
$$
$$=\underbrace{\frac{1}{2}  \lambda  (\lambda + 1)}_{(1)} + 
\underbrace{\frac{1}{2} \left(\frac{n + 3 - 2\lambda}{2}\right)\left(\frac{n + 3 - 2\lambda}{2} - 3\right)}_{(2)} +$$
$$+ \underbrace{\frac{1}{2}\sum\limits_{j=3}^{2 + \frac{n-1}{2} - \lambda} 2  (2 - 2j + 1)}_{(3)}+\underbrace{\frac{1}{2}\sum\limits_{j=3 + \frac{n-1}{2} - \lambda}^{\lambda} 1  (1 - 2j + 1)}_{(4)},$$
where after calculations we obtain:

\begin{enumerate}[(1)]
    \item = $\frac{1}{2}\cdot(\lambda^2 + \lambda);$
    \item = $\frac{1}{8}\cdot(n^2 + 4 \lambda ^ 2 - 4n\lambda - 12\lambda - 9 );$
    \item = $\frac{1}{4}\cdot (-n^2 -4\lambda^2 + 4n\lambda -2n + 4\lambda + 3);$
    \item = $\frac{1}{8}\cdot (n^2 - 4n\lambda + 4n - 4\lambda + 3).$
\end{enumerate}
Finally, putting all the members of the expression together we have:

$$\frac{1}{2}(\lambda^2 + \lambda) + \frac{1}{8}(n^2 + 4 \lambda ^ 2 - 4n\lambda - 12\lambda - 9 ) + $$ $$+\frac{1}{4} (-n^2 -4\lambda^2 + 4n\lambda -2n + 4\lambda + 3) + 
\frac{1}{8} (n^2 - 4n\lambda + 4n - 4\lambda + 3) = \lambda.$$

It is easy to see that this partition holds when $\frac{n + 3}{4} \leqslant \lambda \leqslant \frac{n - 1}{2}$. In addition to the limitations on $\lambda$ there is a limitation on $n$ obtained by substituting $\lambda=\frac{n+3}{4}$ into the expression for partition which gives $n \geqslant 5$ and completes the proof.

\hfill $\square$
\begin{lemma}\label{even_n2}
    If $n\geqslant 10$ and $n$ is even then the partition $(\lambda + 1, \frac{n + 2 - 2\lambda}{2}, 3,  2\times (\frac{n-4}{2} - \lambda), 1\times (\frac{4\lambda - n - 2}{2})$ corresponds to the eigenvalue $\lambda \in \mathbb{N}$, where $\frac{n + 2}{4} \leqslant \lambda \leqslant \frac{n - 4}{2}$.
\end{lemma}

\proof From~(\ref{transp_eigen}) we have: 
$$\lambda_{\left(\lambda + 1, \frac{n + 2 - 2 \lambda}{2}, 3, 2\times(\frac{n-4}{2} - \lambda), 1\times (\frac{4\lambda - n - 2}{2}) \right)}=
$$

$$=\underbrace{\frac{1}{2}  \lambda  (\lambda + 1)}_{(1)} + 
\underbrace{\frac{1}{2} \left(\frac{n + 2 - 2\lambda}{2}\right)\left(\frac{n + 2 - 2\lambda}{2} - 3\right)}_{(2)} + \underbrace{\frac{1}{2}3(3-6+1))}_{(3)}+$$

$$+\underbrace{\frac{1}{2}\sum\limits_{j=4}^{3 + \frac{n-4}{2} - \lambda} 2  (2 - 2j + 1)}_{(4)} + \underbrace{\frac{1}{2}\sum\limits_{j=4 + \frac{n-4}{2} - \lambda}^{\lambda} (1 - 2j + 1)}_{(5)}.$$

Hence, after calculations we get:
\begin{enumerate}[(1)]
    \item = $\frac{1}{2}\cdot(\lambda^2 + \lambda);$
    \item = $\frac{1}{8}\cdot(n^2 + 4 \lambda^2 - 4n\lambda -2n +4\lambda-8);$
    \item = $-3;$
    \item = $\frac{1}{4}\cdot(-n^2-4\lambda^2+4n\lambda+16);$
    \item = $\frac{1}{8}\cdot(n^2-4n\lambda + 2n),$
\end{enumerate}
and finally we obtain:

$$\frac{1}{2}\cdot(\lambda^2 + \lambda) + \frac{1}{8}\cdot(n^2 + 4 \lambda^2 - 4n\lambda -2n +4\lambda-8) - 3 +$$ $$+\frac{1}{4}\cdot(-n^2-4\lambda^2+4n\lambda+16) + \frac{1}{8}\cdot(n^2-4n\lambda + 2n) = \lambda.$$

It is easy to see that this partition holds for any integer $\lambda \in  [\frac{n + 2}{4},\frac{n - 4}{2}]$. Moreover, by substituting $\lambda=\frac{n + 2}{4}$ into the partition we get $n \geqslant 10$. \hfill $\square$

\begin{lemma}\label{even_mod0}
    If $n\geqslant 20$ and $n\equiv 0 \ ({\rm mod} \ 4)$ then the partition $(\frac{n}{4}, \frac{n}{4}, 5, 4, 2 \times (\frac{n-20}{4}), 1)$ corresponds to the eigenvalue $[\frac{n-3}{4}] + 1$.
\end{lemma}

\proof Note that if $n\equiv 0 \ ({\rm mod} \ 4)$ then $[\frac{n-3}{4}] + 1 = \frac{n}{4}$. Substituting the partition into~(\ref{transp_eigen}) we immediately have:

$$\lambda_{\left(\frac{n}{4}, \frac{n}{4}, 5, 4,  2\times (\frac{n - 20}{4}), 1\right)}=
$$

$$=\underbrace{\frac{1}{2} \frac{n}{4}\left(\frac{n}{4} - 2 + 1\right)}_{(1)} + 
\underbrace{\frac{1}{2} \frac{n}{4}\left(\frac{n}{4} - 4 + 1\right)}_{(2)} + 
\underbrace{\frac{1}{2} 5\left(5 - 6 + 1\right)}_{(3)} + 
$$

$$
+ \underbrace{\frac{1}{2} 4\left(4 - 8 + 1\right)}_{(4)}
+ \underbrace{\frac{1}{2}\sum\limits_{j=4}^{j=3+\frac{n-20}{4}}2(2-j+1)}_{(5)} 
+ \underbrace{\frac{1}{2} \left(1-(5+\frac{n-20}{4}) + 1\right)}_{(6)}, 
$$
where each of the terms after calculations is given as follows:
\begin{enumerate}[(1)]
    \item = $\frac{1}{32}\cdot(n^2 - 4n);$
    \item = $\frac{1}{32}\cdot(n^2 - 12n);$
    \item = $0;$
    \item = $-6;$
    \item = $\frac{1}{16}\cdot(-n^2+16n+80);$
    \item = $\frac{1}{4}\cdot(4-n),$
\end{enumerate}
and finally we have:

$$\frac{1}{32}(n^2 - 4n) +
\frac{1}{32}(n^2 - 12n) - 6 + \frac{1}{16}(-n^2+16n+80) + \frac{1}{4}(4-n) = \frac{n}{4}$$
with holding the partition for $n\geqslant 20$. \hfill $\square$

\begin{lemma}\label{odd_mod1}
If $n\geqslant 13$ and  $n\equiv 1 \ ({\rm mod} \ 4)$ then the partition $(\frac{n + 3}{4}, \frac{n + 3}{4}, 4,  2\times (\frac{n - 13}{4}), 1)$ corresponds to the eigenvalue $[\frac{n-3}{4}] + 1$.
\end{lemma}

\proof It is easy to see that if $n\equiv 1 \ ({\rm mod} \ 4)$ then $[\frac{n-3}{4}] + 1 = \frac{n - 1}{4}$. A direct substitution of the partition into~(\ref{transp_eigen}) leads to the following expression:

$$\lambda_{\left(\frac{n + 3}{4}, \frac{n + 3}{4}, 4,  2\times (\frac{n - 13}{4}), 1\right)}=
$$

$$=\underbrace{\frac{1}{2} \left(\frac{n+3}{4}\right)\left(\frac{n+3}{4} - 2 + 1\right)}_{(1)} + 
\underbrace{\frac{1}{2} \left(\frac{n+3}{4}\right)\left(\frac{n+3}{4} - 4 + 1\right)}_{(2)} + 
$$

$$
+\underbrace{\frac{1}{2}\left(4-6+1\right)}_{(3)}
+ \underbrace{\frac{1}{2}\sum\limits_{j=4}^{j + \frac{n-13}{4}} 2 (2 - 2j + 1)}_{(4)} + \underbrace{\frac{1}{2}\left(1-2\frac{n+3}{4} + 1\right)}_{(5)}.
$$
Further calculations give the terms $(1)=\frac{1}{32}\cdot(n^2+2n-3)$; $(2)=\frac{1}{32}\cdot(n^2-6n-27)$; $(3)=-2$; $(4)=-\frac{1}{16}\cdot(n^2-10n-39)$; and $(5)=\frac{1}{4}\cdot(1-n),$ 
whose summation leads to the desired eigenvalue  $\frac{n - 1}{4}$. The only thing left to notice is that the partition holds for any $n\geqslant 13$. \hfill $\square$

\begin{lemma}\label{even_mod2}
If $n \geqslant 22$  and $n\equiv 2 \ ({\rm mod} \ 4)$ then the partition $(\frac{n + 2}{4}, \frac{n - 2}{4}, 5, 4, 2 \times (\frac{n-22}{4}), 1, 1)$ corresponds to the eigenvalue $[\frac{n-3}{4}] + 1$. 
\end{lemma}

\proof It is clear that if $n\equiv 2 \ ({\rm mod} \ 4)$ then  $[\frac{n-3}{4}] + 1=\frac{n - 2}{4}$. By a direct substitution of the partition into~(\ref{transp_eigen}), we immediately have:

$$\lambda_{\left(\frac{n + 2}{4}, \frac{n - 2}{4}, 5, 4,  2\times (\frac{n - 22}{4}), 1, 1\right)}=
$$

$$=\underbrace{\frac{1}{2} \left(\frac{n+2}{4}\right)\left(\frac{n+2}{4} - 2 + 1\right)}_{(1)} + 
\underbrace{\frac{1}{2} \left(\frac{n-2}{4}\right)\left(\frac{n-2}{4} - 4 + 1\right)}_{(2)} + 
$$

$$
\underbrace{\frac{1}{2}5\left(5-6+1\right)}_{(3)}
+ \underbrace{\frac{1}{2}4\left(4-8+1\right)}_{4)}
+ \underbrace{\frac{1}{2}\sum\limits_{j=4}^{j + \frac{n-22}{4}} 2 (2 - 2j + 1)}_{(5)} + 
$$

$$
+\underbrace{\frac{1}{2}\left(1-2\frac{n-2}{4} + 1\right)}_{(6)}
+ \underbrace{\frac{1}{2}\left(1-2\frac{n + 2}{4} + 1\right)}_{(7)},
$$
whose terms after calculations are given by $(1)=\frac{1}{32}\cdot(n^2-4)$; $(2)=\frac{1}{32}\cdot(n^2-16n+28)$; $(3)=0$; $(4)=-6$; $(5)=-\frac{1}{16}\cdot(n^2-20n-44)$; $(6)=\frac{1}{4}\cdot(6-n)$; and $(7)=\frac{1}{4}\cdot(2-n)$, and after summation we have:
$$\frac{1}{32}(n^2 - 4) + \frac{1}{32}(n^2 - 16n + 28) -6 -\frac{1}{16}(n^2-20n-44) + \frac{1}{4}(6-n) + \frac{1}{4}(2-n) = \frac{n - 2}{4}.$$
In this case the partition holds for $n\geqslant 22$.
\hfill $\square$

\begin{lemma}\label{odd_mod3}
If $n\geqslant 19$ and $n\equiv 3 \ ({\rm mod} \ 4)$ then the partition $(\frac{n + 1}{4}, \frac{n + 1}{4}, 5, 3, 2\times (\frac{n - 19}{4}), 1)$ corresponds to the eigenvalue $[\frac{n-3}{4}] + 1$.
\end{lemma}

\proof It is obvious that if $n\equiv 3 \ ({\rm mod} \ 4)$ then $[\frac{n-3}{4}] + 1=\frac{n + 1}{4}$. Using~(\ref{transp_eigen}) we have:

$$\lambda_{\left(\frac{n + 1}{4}, \frac{n + 1}{4}, 5, 3, 2\times (\frac{n - 19}{4}), 1\right)}=
$$

$$=\underbrace{\frac{1}{2} \left(\frac{n+1}{4}\right)\left(\frac{n+1}{4} - 2 + 1\right)}_{(1)} + 
\underbrace{\frac{1}{2} \left(\frac{n+1}{4}\right)\left(\frac{n+1}{4} - 4 + 1\right)}_{(2)} + \underbrace{\frac{1}{2}5\left(5-6+1\right)}_{(3)}+$$

$$
+ \underbrace{\frac{1}{2}3\left(3-8+1\right)}_{(4)}
+ \underbrace{\frac{1}{2}\sum\limits_{j=5}^{5 + \frac{n-19}{4}} 2 (2 - 2j + 1)}_{(5)} + \underbrace{\frac{1}{2}\left(1-2\frac{n+1}{4} + 1\right)}_{(6)}
$$
with reduction to the following terms: $(1)=\frac{1}{32}\cdot(n^2-2n-3)$; $(2)=\frac{1}{32}\cdot(n^2-10n-11)$; $(3)=0$; $(4)=-6$; $(5)=-\frac{1}{16}\cdot(n^2-14n-95)$; and $(6)=\frac{1}{4}\cdot(3-n),$ whose summation gives $\frac{n + 1}{4}$, and the partition holds for $n\geqslant 19$. \hfill $\square$

\section{Discussions and further research}
In this paper we have shown that all integers from the interval $[-\frac{n-4}{2}, \frac{n-4}{2}]$ belong to the spectrum of $T_n$ improving the bounds $[-\frac{n-4}{10}, \frac{n-4}{10}]$ obtained in~\cite{KK22}. However, the full description of the spectrum is still unknown.

Our main goal in future research is to give the most complete description of the spectrum of $T_n$ for large values of $n$. We are sure that the segment $[-\frac{n-4}{2}, \frac{n-4}{2}]$ is not the best estimate and that this estimate can be improved by some larger segment $I(n)$ represented as follows:

$$I(n)=s_0\cup (\cup_{k=1}^{l} s_k ) \cup \{a_1, \dots, a_m\},$$
where $s_0$ and $s_1$ are the segments considered in the proof of Theorem~\ref{segment_0_n2_theorem}. For example, $s_0=[1, \frac{n-3}{4}]$ and $s_1=[\frac{n+3}{4}, \frac{n-1}{2}]$ if $n$ is odd. Hence, for segments $\cup_{k=1}^ls_k$ there should be a more general description of partitions like this: $$p(n, k, \lambda)=(n_1(n, k, \lambda), n_2(n, k, \lambda), \dots, n_r(n, k, \lambda)).$$ It would clearly follow from such a representation that all eigenvalues from the set $\cup_{k=1}^ls_k$ belong to the spectrum of the Transposition graph $T_n$. We distinguish the segment $s_0$ among the other segments since its partition has the first part depending on $\lambda$ and $n$ though we assume the existence of a representation $p(n, k, \lambda)$ whose the first part does not depend on $n$, i.~e. $n_1(n, k,\lambda)=n_1(k, \lambda)$.

\section*{Acknowledgements}
The study by Elena~V.~Konstantinova was performed according to the Government research assignment for IM SB RAS, project FWNF-2022-0017. The work of Artem Kravchuk was supported by the Mathematical Center in Akademgorodok, under agreement No. 075-15-2022-281 with the Ministry of Science and High Education of the Russian Federation.

\end{document}